\documentclass[preprint,12pt]{elsarticle}

\usepackage{graphicx}
\usepackage{amssymb}

\usepackage{lineno}

\usepackage{comment}
\usepackage{ragged2e}
\usepackage{xcolor}
\usepackage{colortbl}
\usepackage{soul}
\usepackage{forest}
\usepackage{algorithm}
\usepackage{algorithmic}
\usepackage{amsmath}
\usepackage{mathtools}
\usepackage{bbm}
\usepackage{dirtytalk}
\useforestlibrary{edges}
\usepackage{mathtools}
\usepackage{tabularx}
\usepackage{multirow,tabularx}
\usepackage{caption}
\usepackage{color}
\usepackage{multirow}
\usepackage{array}
\usepackage{tabularx}
\usepackage{float}
\usepackage{setspace}
\usepackage{lipsum}
\usepackage{hyperref}
\usepackage{graphicx}
\usepackage{subfig}
\usepackage{booktabs}
\usepackage[export]{adjustbox}

\journal{Journal of Computers $\&$ Operations Research}

\begin{document}

\begin{frontmatter}






\title{Faster Maximum Feasible Subsystem Solutions for Dense Constraint Matrices}

\author{Fereshteh Fakhar Firouzeh\corref{mycorrespondingauthor}}
\cortext[mycorrespondingauthor]{Corresponding author}
\ead{behnazfakharfirouzeh@cmail.carleton.ca}

\author{John W. Chinneck}

\author{Sreeraman Rajan}

\address{Systems and Computer Engineering, Carleton University, Ottawa, Ontario, Canada}

\begin{abstract}
Finding the largest cardinality feasible subset of an infeasible set of linear constraints is the \textit{Maximum Feasible Subsystem} problem (MAX FS). Solving this problem is crucial in a wide range of applications such as machine learning and compressive sensing. Although MAX FS is NP-hard, useful heuristic algorithms exist, but these can be slow for large problems. We extend the existing heuristics for the case of dense constraint matrices to greatly increase their speed while preserving or improving solution quality. We test the extended algorithms on two applications that have dense constraint matrices: binary classification, and sparse recovery in compressive sensing. In both cases, speed is greatly increased with no loss of accuracy. 

\end{abstract}

\begin{keyword}
Maximum Feasible Subsystem \sep Linear Programming \sep Classification \sep Compressive Sensing


\end{keyword}

\end{frontmatter}


\section{Introduction}
Finding the maximum cardinality feasible subsystem of an infeasible set of linear constraints is known as the \textit{Maximum Feasible Subsystem} problem (MAX FS) \cite{One}. Its complement is the \textit{Minimum Unsatisfied Linear Relation} problem (MIN ULR) \cite{amaldi1994finding} of finding the minimum number of constraints to remove from an infeasible set such that the remaining constraints have a feasible solution. Moreover, all infeasible systems have atleast one \textit{Irreducible Infeasible Subsets} (IISs) in their original set of constraints. An IIS is a set of constraints that is infeasible, but is rendered feasible if any constraint is removed. If all of the IISs in the model are enumerated, then a MIN ULR solution is found as a minimum cover of the IISs; this is known as the \textit{minimum IIS cover} problem \cite{chinneck1996effective}. In this paper, MAX FS, MIN ULR, and MIN IIS COVER are the same problem and the terms will be used interchangeably.

MAX FS is NP-hard \cite{Seven,Eight,Nine} so exact solution algorithms are only available for relatively small instances, but solutions are still needed in a wide variety of applications including machine learning \cite{Three}, misclassification minimization \cite{Four}, training of neural networks \cite{amaldi1994finding}, telecommunications \cite{Five}, computational biology \cite{Six}, and compressive sensing (CS) \cite{firouzeh2020maximum}. To be useful in practice, MAX FS heuristics should obtain high quality solutions (large feasible subsets, or equivalently, small MIN ULR subsets) quickly.

The MAX FS problem can be formulated for exact solution in various ways, including as a mixed-integer linear program \cite{chinneck2008feasibility}, a linear program with equilibrium constraints (LPEC) \cite{amaldi2003maximum,mangasarian1994misclassification,mangasarian1996machine,bennett1997parametric}, or as a type of set-covering problem \cite{parker1995set,parker1996finding,pfetsch2003maximum,pfetsch2008branch}. As shown in chapter 7 of \cite{chinneck2008feasibility}, there are practical difficulties in using any of these exact methods to solve this NP-hard problem (i.e. they can be used for a very small MAX FS problem). For this reason, heuristic solutions are required. 

Chinneck \cite{chinneck1996effective,chinneck2001fast} developed a set of effective linear programming (LP)-based heuristics for solving the MAX FS problem, and showed experimentally that they give better results than  methods such as MISMIN, a state-of-the-art LPEC solver at the time \cite{chinneck2001fast}. Variants of Chinneck's methods have been developed more recently for finding sparse solutions to underdetermined systems of linear equations in unbounded variables \cite{firouzeh2020maximum} for sparse recovery in compressive sensing.

This paper extends Chinneck's algorithms to provide major improvements in solution speed for the case of dense constraint matrices, based on the observation that dense constraint matrices provide multiple similar candidates for removal during the MAX FS solution process. The new extensions are evaluated on the following two applications that have dense constraint matrices: binary classification, and sparse recovery in compressive sensing.

Chinneck's algorithms remove constraints from the original set one by one until what remains is the heuristic MAX FS solution; the removed constraints constitute the heuristic MIN ULR solution. An inner loop assesses candidate constraints, and an outer loop removes the single constraint chosen by the inner loop.  Two extensions that greatly improve solution speed when a dense constraint matrix provides multiple similar candidates are proposed in this paper. Extension 1 amalgamates the two loops into a single one that both assesses the candidates and removes \textit{multiple} constraints at each iteration. Extension 2 provides an early exit under certain conditions. When tested on binary classification and compressive sensing sparse recovery, the extended heuristics greatly reduce solution time while preserving or improving solution quality.  

The paper is structured as follows. Section \ref{survey} reviews the existing MAX FS solution heuristics. The new extensions are developed in Section \ref{Novel}. The experimental results for binary classification are given in Section \ref{Classification} and for compressive sensing sparse recovery in Section \ref{USLE}. Section \ref{Conc} concludes the paper and outlines future work.

\section{LP-based MAX FS Solution Heuristics} \label{survey}

Chinneck's original polynomial time heuristics \cite{chinneck1996effective} are briefly summarized here. All variants first elasticize the infeasible set of linear row constraints (and possibly variable bounds) by adding nonnegative \textit{elastic variables} as shown in Table \ref{ElasticVersions}. The elastic variables measure the extent to which each constraint is violated, so a linear program (LP) minimizing their sum finds the minimum total violation for an infeasible set of constraints and bounds. There are two types of elastic models: the \textit{standard elastic LP} elasticizes only the row constraints, and the \textit{full elastic LP} elasticizes both the row constraints and the column bounds. The objective for a full elastic LP is: {$min\, Z = \sum_{i} (e_i^+ + e_i^-)+\sum_{j} (e_j^+ + e_j^-)$} for constraints $i = 1...m$ and bounded variables $j= 1...n$; the second term is omitted for a standard elastic LP.  $Z=0$ indicates a feasible system. 

\begin{table}[hbt!]
\caption{Elasticizing constraints by adding non-negative elastic variables.}
\footnotesize
\centering

\centering
\begin{tabular}{||c||l|l|l||}
\hline
\label{ElasticVersions}

\textbf{Type} & \textbf{Nonelastic} & \textbf{Standard elastic} & \textbf{Full elastic} \\ \hline\hline

\multirow{3}{*}{\textbf{Row Cons.}}& $\sum_j a_{ij} x_j \geq {b}_i$& $\sum_j a_{ij} x_j + e_i \geq {b}_i $& $\sum_ja_{ij} x_j + e_i \geq {b}_i$\\ \cline{2-4} 
  
& $\sum_j a_{ij} x_j \leq {b}_i$&  $\sum_j a_{ij} x_j - e_i \leq {b}_i  $& $\sum_j a_{ij} x_j - e_i \leq {b}_i $\\ \cline{2-4} 

& $\sum_j a_{ij} x_j = {b}_i$ & $a_{ij} x_j + e_i^+ - e_i^- = {b}_i$ &  $\sum_ja_{ij} x_j + e_i^+ - e_i^- = {b}_i$\\ \hline

\multirow{2}{*}{\textbf{Var. Bnds}} &$x_j \geq l_j$&$x_j \geq l_j$&$x_j+e_j^+ \geq l_j$\\ \cline{2-4} 

&$x_j \leq u_j$&$x_j \leq u_j$&$x_j - e_j^-\leq u_j$\\ \hline

\end{tabular}
\end{table}

As shown in Fig. 1, the General MAX FS algorithm of \cite{chinneck1996effective} has an inner loop (line 10) and an outer loop (line 5). The inner loop tests each member of a set of \textit{candidate constraints} in \textit{CandidateSet} by examining their effect on $Z$ when temporarily removed from the model. The candidate constraint that most reduces $Z$ when removed is selected as the \textit{Winner} candidate, and is subsequently removed from the model and added to \textit{MINULR} in the outer loop. This process continues until $Z=0$ (feasibility is reached). The constraints remaining in the model constitute the heuristic MAX FS solution; the removed constraints in \textit{MINULR} constitute the heuristic MIN ULR solution. Algorithm variants differ in how the list of candidate constraints is constructed.

The heuristic relies on \textit{Observation 3} from \cite{chinneck1996effective}, which notes that a constraint that appears in more than one IIS will reduce $Z$ more than a constraint that appears in a single IIS because it will eliminate multiple IISs upon removal. This means that a large drop in $Z$ when a constraint is removed is a useful sign that the constraint is part of the MINULR set. The larger the drop in $Z$, the more IISs the constraint likely covers.
 
There are efficiency enhancements in terms of algorithm processing time. If the  \textit{CandidateSet} has a single element (line 6), then that element is added to \textit{MINULR} and the algorithm exits: there is no need to solve the LP. As the $Winner$ constraint is updated, so is the list of candidates for the next round, $NextCandidateSet$. All versions of the algorithm solve multiple LPs, but each LP is identical to the last one except for two constraints, so simplex advanced starts to speed the process considerably.

Algorithm 1 to 3 of the base algorithm differ in how the list of candidate constraints is constructed, as described next.

\subsection{Algorithm 1}
Candidate constraints are those to which the elastic objective function is sensitive (nonzero dual price). Constraints not in this list do not affect $Z$ when dropped, so they are not candidates (\textit{Observation 4} from \cite{chinneck1996effective}). This is the longest possible candidate list, but its length can be limited as follows. Sort the constraints in descending order by their absolute dual cost and take only the first $k$ elements of the sorted list. This is referred to as Algorithm 1$(k)$. For clarity in the sequel, unlimited lists are indicated by $k=\infty$.

\begin{algorithm}
\captionof{figure}{General MAX FS Algorithm \cite{chinneck1996effective}}
\begin{algorithmic} [t]
\footnotesize
	\sffamily
\STATE $\textbf{1:} \textit{ MINULR} \leftarrow \emptyset$ \\%
\textbf{2:}\text{ Set up elastic LP.}\\
\textbf{3:}\text{ Solve elastic LP.}\\
\textbf{4:}\text{ Construct }\textit{CandidateSet.}\\
\textbf{5:}\textbf{ do} \text{[\textit{outer loop}]:}  \\
\textbf{6:}\hspace{1em} \textbf{if } $\mid \textit{CandidateSet} \mid $ \text{= 1 } \textbf{ then}\\
\textbf{7:}\hspace{2em}\text{ Add the single candidate to} \textit{MINULR} \text{and exit}.\\
\textbf{8:}\hspace{1em}\textbf{ end if}\\
\textbf{9:}\hspace{1em}\textit{ WinnerZ} $\leftarrow \infty.$ \\%
\textbf{10:}$\hspace{1em}\textbf{for} \text{ each candidate $j$ in } \textit{CandidateSet} \textbf{ do } \text{[\textit{inner loop}]:}$\\
\textbf{11:}$\hspace{2em}\text{Delete candidate $j$.}$\\
\textbf{12:}$\hspace{2em}\text{Solve elastic LP.}$\\
\textbf{13:}$\hspace{2em}\textbf{if } Z = 0 \textbf{ then}$\\
\textbf{14:}$\hspace{3em} \text{Add candidate $j$ to \textit{MINULR} and exit.}$\\
\textbf{15:}$\hspace{2em}\textbf{end if}$\\
\textbf{16:}$\hspace{2em}\textbf{if } Z < \textit{WinnerZ}  \textbf{ then}$\\
\textbf{17:}$\hspace{3em}\textit{Winner} \leftarrow \text{candidate } j.$\\
\textbf{18:}$\hspace{3em}\textit{WinnerZ} \leftarrow Z.$\\
\textbf{19:}$\hspace{3em}\text{Construct } \textit{NextCandidateSet}.$\\
\textbf{20:}$\hspace{2em}\textbf{end if}$\\
\textbf{21:}$\hspace{2em}\text{Reinstate candidate $j$.}$\\
\textbf{22:}$\hspace{1em}\textbf{end for}$\\
\textbf{23:}$\hspace{1em}\text{Add \textit{Winner} to \textit{MINULR}.}$ \\%
\textbf{24:}$\hspace{1em} \text{Delete \textit{Winner} permanently.}$\\
\textbf{25:}$\hspace{1em} \textit{CandidateSet} \leftarrow \textit{NextCandidateSet}.$\\
\textbf{26:}$ \textbf{ end do}$
\STATE \textbf{OUTPUT}: \textit{MINULR} is a heuristic MIN ULR solution. 
\label{genMAXFS}
\end{algorithmic}
\end{algorithm}

\subsection{Algorithm 2} Algorithm 2 uses a different criterion for constructing the \textit{CandidateSet} based on \textit{Observation 5} in \cite{chinneck2001fast}: a good predictor of the magnitude of the drop in $Z$ obtained by deleting a violated constraint $j$ is given by Eqn. (\ref{Maxprod}): 

\begin{equation}
\textit{Product}: \textit{$e_j$} \times |\textit{constraint j dual price}|
\label{Maxprod}
\end{equation}
For equality constraints, the larger of the two elastic variables is used in Eqn. (\ref{Maxprod}). Algorithm 2 can also limit the length of the candidate list. First sort the candidate constraints in descending order by product magnitude, and then add only the top $k$ constraints to \textit{CandidateSet}; this is Algorithm 2($k$). List length $k=1$ is particularly useful in some applications, where it is referred to as the \textit{maximum product} algorithm.

\subsection{Algorithm 3}
Algorithm 3 uses \textit{Observation 6} \cite{chinneck2001fast}: for satisfied constraints, a good predictor of the \textit{relative} magnitude of any drop in $Z$ that may be obtained by deleting the constraint is given by $|\textit{constraint dual price}|$. Algorithm 3  uses 2 lists: (\textit{i}) for violated constraints, the list of product magnitudes (Eqn. (\ref{Maxprod})), and for satisfied constraints (\textit{ii}) the list of absolute constraint sensitivities. Both lists are independently sorted in descending order, and the top $k$ elements are taken from each list, resulting in $2k$ candidates in total. This is referred to as Algorithm 3(\textit{k}).

\subsection{Finding Sparse Solutions to Linear Systems} \label{findingsparsesolutions}

A sparse solution is one having very few nonzeros. Algorithms 1-3 can be used to find sparse solutions to general systems of linear constraints as follows. Given the system $\mathbf{Ax \{\leq, =, \geq\} b}$, add the variable zeroing constraints $\mathbf{x}=\mathbf{0}$. Now apply any of Algorithms 1-3 with candidates chosen only from among the constraints $\mathbf{x}=\mathbf{0}$. The MAX FS solution provides the sparsest solution by allowing the smallest number of nonzero variables. Finding the sparsest solution to an underdetermined system of linear equalities in unbounded variables is of particular interest in a number of applications. We focus on this problem here. 

Chinneck's algorithms have been adapted for underdetermined systems of linear equations in unbounded variables. Jokar and Pfetsch \cite{jokar2008exact} modified Chinneck's algorithm as follows. Given the $m\times n$ linear system $\mathbf{Ax=b}$ in unbounded variables, all variables $x_j$ are first replaced by $u_j-v_j$ where both $u_j$ and $v_j$ are nonnegative, and the following LP is constructed:
\begin{equation}
min\, Z = \sum_{j} (u_j + v_j) \quad s.t. \quad \mathbf{A (u-v)=b, \quad u \geq0, \hspace{0.2em} v \geq 0}
\label{JP}
\end{equation}
This LP has $m$ constraints in $2n$ nonnegative variables. The objective function attempts to drive all variables to zero in a manner similar to driving elastic variables to zero in Algorithms 1-3. 

The basic algorithm follows the logic in Fig. 1, except that candidates are variables and not constraints. Additional differences from Fig. 1:

\begin{itemize}

\item Candidate variables are those not in the MIN ULR set and having a nonzero value in the current LP solution. Candidates are sorted in decreasing order by the magnitude of $|u_j - v_j|$. 
  
\item Candidate variables are tested in the inner loop by temporarily removing them from the objective function by setting the objective coefficients of the associated $u_j$, $v_j$ pair to zero and re-solving the LP.

\item A candidate variable $x_j$ is moved into the MIN ULR set by resetting the objective function coefficients of $u_j$ and $v_j$ to zero permanently so that $x_j > 0$ no longer affects $Z$.
  
\item At exit, the variables in the MIN ULR set are those that can take nonzero values in the sparse solution.

\item Extraneous variables are sometimes included in the MIN ULR solution. A simple post-processing step may remove them. First construct $\mathbf{A}_{mn}$ by eliminating all columns from $\mathbf{A}$ that are not in the MIN ULR solution. For each variable in the MIN ULR set in turn, force it to zero and test the feasibility of $\mathbf{A}_{mn}\mathbf{x=b}$; if feasible then the tested variable is permanently removed from the MIN ULR set.
  
\end{itemize}

The well-known $Basis Pursuit$ (BP) algorithm \cite{chen2001atomic} solves system (\ref{JP}) a single time. This works well if the solution is very sparse, but when this is not the case, Fig.1 with the modifications above is more effective.

Several other recent variants have been shown to be very effective for this problem \cite{ChinneckMAXFS,firouzeh2020maximum}:  
\begin{itemize}

\item \textbf{Method B} uses system (\ref{JP}). The objective  coefficients of the $u_j$ and $v_j$ associated with variable $x_j$ moved into the MIN ULR set are permanently reset to 0.1 instead of zero \cite{ChinneckMAXFS}. Since $Z$ never reaches 0, exit is instead triggered when $\mid\textit{CandidateSet}\mid = 0$.   
    
\item \textbf{Method M} combines Method B with \textit{Basis Pursuit} (BP) \cite{chen2001atomic}, based on the observation that if BP fails, it returns a solution that is typically of size $m$ or nearly so. Thus,  the BP result (the set of nonzeros in the first solution of (\ref{JP})) is taken when the number of nonzeros is less than $m-3$. Method B is applied when the number of nonzeros in that first solution is larger than $m-3$.
    
 \item \textbf{Method C} has explicit elastic variable zeroing constraints, resulting in the following LP:
\begin{equation}
\min Z = \sum_{j} (e_j^+ + e_j^-) \, s.t. 
\!
\left [
\begin{array}{c c c }
\mathbf{A} & \mathbf{0}_{m \times n} & \mathbf{0}_{m \times n}\\
\mathbf{I}& \mathbf{I} & \mathbf{-I} \\
\end{array}
\right ] \!
\left [
\begin{array}{c}
\mathbf{x} \\
\mathbf{e^+}\\
\mathbf{e^-}\\
\end{array}
\right ] =
\left [
\begin{array}{c }
\mathbf{y}_{m \times 1} \\
\mathbf{0}_{n \times 1}\\ 
\end{array}
\right ] \, 
\label{eq:kron_measurement}
\end{equation}
where $e_j^+$ and $e_j^-$ are nonnegative, $\mathbf{A}$ is $m \times n$ and $\mathbf{I}$ is $n \times n$. This LP has $m+n$ constraints in $3n$ variables. There are two lists of candidates: (1) a list based on the magnitude of nonzero variables not yet assigned to the MIN ULR set, and (2) a list based on the dual prices of the variable zeroing constraints $x_j+e_j^+-e_j^-=0$. The selected variable is added to the MIN ULR set by zeroing the objective  coefficients of the corresponding elastic variables.

\end{itemize}

\section{Extensions to Increase Solution Speed} \label{Novel}
The existing algorithms are effective, but there are opportunities for improving their speed. Limiting the length of the candidate list is helpful, but there are still two bottlenecks: evaluating multiple candidates in the inner loop, and removing just a single constraint in the outer loop. New Extension 1 eliminates the inner loop entirely and removes multiple candidate constraints simultaneously in the single remaining  loop. The trick is in identifying which constraints to remove in each iteration. New Extension 2 provides an early exit under certain conditions.

Each of the extensions can be combined with any of the existing list construction algorithms to create new combinations, e.g. combining  Algorithm 2($\infty$) and Extension 1  creates Algorithm 2($\infty$)E1.


\subsection{Extension 1(E1)} \label{E1}

There are three ways to assemble the candidate list in Algorithms 1-3. Regardless of the type of list, previous research shows that the order in the sorted list of candidates is important: those earlier in the list are much more likely to be part of the MIN ULR than those later in the list. Thus, if we remove multiple constraints at each iteration of a single loop, it makes sense to select them from the start of the sorted list. This argument is especially potent when the constraint matrix is dense. In a dense constraint matrix, all constraints involve all of the variables, so all of the constraints are similar in structure. A consequence is that many of the constraints have similar impacts on $Z$. For example, in the classification application, where the constraints are derived from the data points, data points that are close together generate similar constraints. These similar constraints have similar ranking scores in the candidate list, and hence have similar effects on $Z$ when removed from the LP, leading to this \textbf{Observation}: \textit{When candidates have similar ranking measures as the top-ranked element in the candidate list, they are likely to appear in the final MIN ULR set}. 

This observation allows both the elimination of the inner loop to test individual candidate constraints and the simultaneous removal of multiple candidates in the single remaining loop. The question is how to determine which ranking scores are ``similar" to the top score. We sort the candidates by descending ranking score and apply standard methods for detecting abrupt changes in the mean of subsets of the set \cite{killick2012optimal}. We select for removal all candidates up to the first identified abrupt change.

%
%
    

\subsection{Extension 2(E2)}\label{E4}
Extension 2 provides an early exit. If the number of candidates is less than or equal to $\ell$, then permanently remove all candidates and exit the algorithm. This is designated E2($\ell$). 

\section {Application: Binary Classification}\label{Classification}
 Amaldi \cite{amaldi1994finding}, Parker \cite{parker1995set}, Chinneck \cite{chinneck2009tailoring,chinneck2012integrated,chinneck2001fast}, and Silva \cite{silva2017optimization} have shown how the binary classification problem can be transformed into a MAX FS problem: 
 \begin{itemize}
     \item [] \textbf{Given} a training set of $I$ data points $(i=$1…$I)$ in $J$ dimensions $(j=$1…$J)$ where the class of each point is known (Type 0 or 1).
     \item [] \textbf{Define} a constraint for each data point:
        \begin{itemize}
            \item For points of Type 0: $\sum_j d_{ij}w_j \leq w_0 - \epsilon$
            \item For points of type 1: $\sum_j d_{ij}w_j \geq w_0 + \epsilon$
        \end{itemize}
where $d_{ij}$ indicates the value of attribute $j$ for point $i$. The $d_{ij}$ are known constants, $\epsilon$ is a small positive constant and the variables $w_j$ are unrestricted. 
 \end{itemize}
 
Most datasets cannot be completely separated by a single hyperplane, so the initial system is infeasible. After the MAX FS solution, the final LP is feasible, and its solution returns the separating hyperplane $\sum_j d_{ij}w_j = w_0$ which correctly classifies all of the data points in the MAX FS subset while those in the MIN ULR subset are incorrectly classified. This is heuristically the most accurate hyperplane.

\subsection{Binary Classification Algorithm} \label{ClassProAlg}

Prior work shows that Algorithm 2 provides a major speed-up over Algorithm 1 for binary classification \cite{chinneck2001fast}, with a tiny loss in average accuracy. Violated constraints correspond to misclassified points, so the product ranking score (Eqn. \ref{Maxprod}) used in Algorithm 2 estimates the drop in $Z$ relatively accurately. A standard elastic program is used because the variables are unrestricted. Since all attribute values are specified for all data points, the constraint matrix, $\mathbf{D}$, is dense; hence, Extension 1 can be applied. The revised algorithm is designated as Algorithm 2E1, and is summarized in Fig. \ref{Alg2E1}.

\begin{algorithm}
\captionof{figure}{Algorithm 2E1}
\begin{algorithmic} [t]
\footnotesize
	\sffamily
\STATE$\textit{MINULR} \leftarrow \emptyset$ \\%
\text{Set up the standard elastic LP.}\\
\textbf{do:}\\
$\hspace{1em}$\text{Solve LP.}\\
$\hspace{1em}$\textbf{if }$ Z = 0 \textbf{ then}$\\
$\hspace{2em}$\text{Exit.}\\
$\hspace{1em}$\textbf{end if}\\
$\hspace{1em}$\text{Construct }\textit{CandidateSet} using Eq. \ref{Maxprod}.\\
\hspace{1em}\textbf{if } $\mid \textit{CandidateSet} \mid $ \text{= 1 } \textbf{ then}\\
\hspace{2em}\text{ Add the single candidate to} \textit{MINULR} \text{and exit}.\\
\hspace{1em}\textbf{end if}\\
\hspace{1em}Sort candidates from largest to smallest product.\\
\hspace{1em}Select top $p$ candidates for removal based on abrupt change in product score. \\%
\hspace{1em}Add the $p$ selected candidates to \textit{MINULR} and  remove them from the LP. \\%
$\textbf{ end do}$
\STATE \textbf{OUTPUT}: \textit{MINULR} is a heuristic MIN ULR solution. 
\label{Alg2E1}
\end{algorithmic}
\end{algorithm}

\begin{figure}[htp]
\centering
\includegraphics[width=10cm, height = 6.5cm]{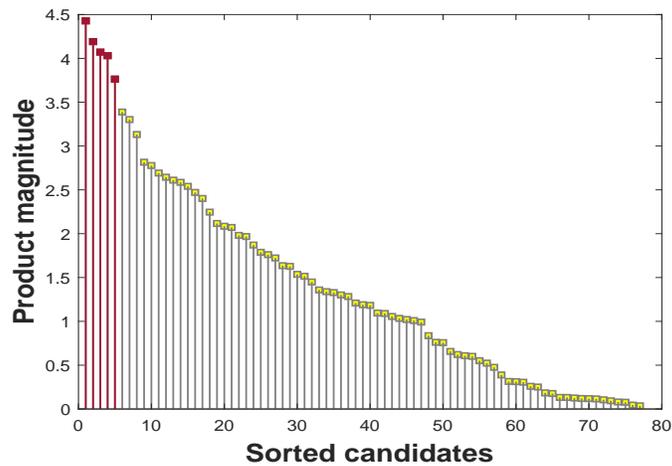}
\caption{\small{First Sorted Candidate List for ``Ozone Level Detection" Dataset.}}
\label{ltr}
\end{figure}

To illustrate the workings of Extension 1, Fig. \ref{ltr} shows the first sorted candidate list for the ``Ozone Level Detection" dataset. The first sharp drop-off occurs after the fifth candidate, so the first five constraints/points are removed in this case. 

\subsection{Comparators} 
We compare Algorithm 2E1 with two existing algorithms: (i) original Algorithm 2($\infty$) and (ii) original Algorithm 2(1). Both have been shown to provide high accuracies in binary classification \cite{chinneck2001fast}.

\subsection{Evaluation Metric}
Accuracy is used to evaluate the quality of the solutions:

 \begin{equation}
     \textit{Accuracy} = \frac{\textit{Total $\#$ correctly classified instances}}{\textit{Total $\#$  instances}}
     \label{accuracyEQ}
 \end{equation}
For comparison with previous results \cite{chinneck2001fast}, the training set is identical to the entire dataset.

\subsection{Data Sets} 
Twenty binary classification problems are derived from datasets in the UCI Repository of Machine Learning Databases \cite{Dua:2019}. Table \ref{ClassiDataset} summarizes their characteristics.    

         \begin{table}[h]
         \caption{Classification Datasets}
         \centering
         \footnotesize
         \label{ClassiDataset}
    \begin{tabular}{|l||c|c|c|}
    \hline
    \multicolumn{1}{|c||}{\textbf{Dataset}} & \textbf{Instances} & \textbf{Features}  \\ \hline\hline
    Breast Cancer Wisconsin Original       & 683                & 9\\ \hline
    Cardiotocography (Type 1 vs. others)                       & 2126               & 34 \\ \hline
    
    Cardiotocography (Type 2 vs. others)                       & 2126               & 34 \\ \hline
    
    Car Evaluation (Type 1 vs. others)                         & 1728               & 6 \\ \hline
    
    Car Evaluation (Type 2 vs. others)                         & 1728               & 6 \\ \hline
    
    Car Evaluation (Type 3 vs. others)                         & 1728               & 6 \\ \hline
    
    Car Evaluation (Type 4 vs. others)                         & 1728               & 6 \\ \hline
    
    Climate Model Simulation Crashes       & 540                & 18                \\ \hline
    Credit Approval                        & 653                & 15                \\ \hline
    Diabetic Retinopathy Debrecen          & 1151               & 19                \\ \hline
    Glass Identification                   & 214                & 9                 \\ \hline
    Heart Disease                & 297                & 13                \\ \hline
    Hepatitis Domain                       & 112                & 18                \\ \hline
    Ionosphere                             & 351                & 33                \\ \hline
    BUPA Liver Disorders                   & 341                & 6                 \\ \hline
    Ozone Level Detection                   & 1848               & 102                \\ \hline
    Parkinson's Data Set                    & 195                & 22                \\ \hline
    Pima Indians Diabetes                  & 768                & 8                 \\ \hline
    Thyroid Gland Data                     & 215                & 5                 \\ \hline
    Vowel Recognition Data                 & 990                & 11                \\ \hline
    \end{tabular}
    \end{table}

\subsection{Software and Hardware} \label{hwsw}
All algorithms are implemented in Matlab version $2020$.The linear programming solver is MOSEK via the MOSEK Optimization Toolbox for Matlab version $8.1.0.56$ \cite{MOSEK}. Abrupt changes in the ranking scores in the candidate list are identified using the Matlab \textit{ischange}() function. The computations are carried out on a $3.40$ GHz Intel core $i7$ machine with $16.0$ \textit{GB} RAM, running Windows $10$. 

\subsection{Results and Discussion}

The results in this section show the higher accuracy and greater speed of Algorithm 2E1 vs. original Algorithms 2($\infty$) and 2(1) for binary classification. Table \ref{classificationACC} reports algorithm accuracies with the best results bolded. Table \ref{classificationAnalysisAcc} shows the difference from the best accuracy on each data set. The last three rows give additional metrics: the largest difference from the best accuracy, the average difference from the best accuracy, and number of best accuracies obtained by each algorithm. 

Algorithm 2E1 provides better total accuracy. It has the smallest average difference from best accuracy, and provides the best accuracy for $12$ of the $20$ data sets while it was $8$ for the two existing algorithms. The largest difference from the best accuracy for Algorithm 2E1 is $-1.38$, which is smaller than those for the other two ($-4.28$ for both). 

\begin{table}[hbt!]
\centering
\scriptsize
\caption{Algorithm Accuracies for Classification Datasets.}
\label{classificationACC}
\begin{tabular}{|l|c|c|c|}
\hline
\multicolumn{1}{|c|}{\multirow{2}{*}{\textbf{Dataset}}} & \multicolumn{3}{c|}{\textbf{Accuracy}}                                                                                                                        \\ \cline{2-4} 
\multicolumn{1}{|c|}{}                                  & \textbf{Algorithm 2($\infty$)} & \textbf{\begin{tabular}[c]{@{}c@{}}Algorithm 2(1)\end{tabular}} & \textbf{\begin{tabular}[c]{@{}c@{}}Algorithm 2E1\end{tabular}} \\ \hline
Breast Cancer Wisconsin Original                        & 98.10               & 98.10                                                              & \textbf{98.39}                                                     \\ \hline
Cardiotocography (Type 1 vs.others)                     & 99.72               & 99.72                                                              & \textbf{99.86}                                                     \\ \hline
Cardiotocography (Type 2 vs.others)                     & 99.67               & 99.58                                                              & \textbf{99.81}                                                     \\ \hline
Car Evaluation (Type 1 vs.others)                       & 78.47               & 78.47                                                              & \textbf{82.75}                                                     \\ \hline
Car Evaluation (Type 2 vs.others)                       & \textbf{96.01}      & \textbf{96.01}                                                     & \textbf{96.01}                                                     \\ \hline
Car Evaluation (Type 3 vs.others)                       & 87.79               & 88.42                                                              & \textbf{89.01}                                                     \\ \hline
Car Evaluation (Type 4 vs.others)                       & \textbf{98.15}      & 98.09                                                              & 98.09                                                              \\ \hline
Climate Model Simulation Crashes                        & \textbf{98.89}      & 98.70                                                              & 98.52                                                              \\ \hline
Credit Approval                                         & \textbf{88.05}      & 86.68                                                              & 86.68                                                              \\ \hline
Diabetic Retinopathy  Debrecen                          & 79.06               & \textbf{79.41}                                                     & 78.62                                                              \\ \hline
Glass Identification                                    & 84.11               & \textbf{84.56}                                                     & 83.18                                                              \\ \hline
Heart Disease Cleveland                                 & \textbf{89.56}      & \textbf{89.56}                                                     & \textbf{89.56}                                                     \\ \hline
Hepatitis Domain                                        & 92.86               & \textbf{96.43}                                                     & \textbf{96.43}                                                     \\ \hline
Ionosphere                                              & 98.01               & \textbf{98.01}                                                     & 97.73                                                              \\ \hline
BUPA Liver  Disorders                                   & 75.37               & \textbf{75.66}                                                     & 74.78                                                              \\ \hline
Ozon Level  Detection                                   & 98.86               & 99.13                                                              & \textbf{99.19}                                                     \\ \hline
Parkinsons Data Set                                     & 94.36               & 94.87                                                              & \textbf{96.41}                                                     \\ \hline
Pima Indians  Diabetes                                  & \textbf{80.47}      & 80.08                                                              & 79.69                                                              \\ \hline
Thyroid Gland Data                                      & \textbf{94.42}      & \textbf{94.42}                                                     & \textbf{94.42}                                                     \\ \hline
Vowel Recognition Data                                  & 98.38               & 98.38                                                              & \textbf{98.48}                                                     \\ \hline
\end{tabular}
\end{table}


\begin{table}[hbt!]
\centering
\scriptsize
\caption{Differences from Best Accuracy.}
\label{classificationAnalysisAcc}
\begin{tabular}{|l|c|c|c|}
\hline
\multicolumn{1}{|c|}{\textbf{Dataset}}      & \textbf{Algorithm 2($\infty$)} & \textbf{Algorithm 2(1)} & \textbf{Algorithm 2E1} \\ \hline
Breast Cancer Wisconsin Original            & -0.29               & -0.29                   & \textbf{0}             \\ \hline
Cardiotocography (Type 1 vs.others)         & -0.14               & -0.14                   & \textbf{0}             \\ \hline
Cardiotocography (Type 2 vs.others)         & -0.14               & -0.23                   & \textbf{0}             \\ \hline
Car Evaluation (Type 1 vs.others)           & -4.28               & -4.28                   & \textbf{0}             \\ \hline
Car Evaluation (Type 2 vs.others)           & \textbf{0}          & \textbf{0}              & \textbf{0}             \\ \hline
Car Evaluation (Type 3 vs.others)           & -1.22               & -0.59                   & \textbf{0}             \\ \hline
Car Evaluation (Type 4 vs.others)           & \textbf{0}          & -0.06                   & -0.06                  \\ \hline
Climate Model Simulation Crashes            & \textbf{0}          & -0.19                   & -0.18                  \\ \hline
Credit Approval                             & \textbf{0}          & -1.37                   & -1.37                  \\ \hline
Diabetic Retinopathy  Debrecen              & -0.35               & \textbf{0}              & -0.79         \\ \hline
Glass Identification                        & -0.45               & \textbf{0}              & -1.38                  \\ \hline
Heart Disease Cleveland                     & \textbf{0}          & \textbf{0}              & \textbf{0}             \\ \hline
Hepatitis Domain                            & -3.57               & \textbf{0}              & \textbf{0}             \\ \hline
Ionosphere                                  & \textbf{0}          & \textbf{0}              & -0.28                  \\ \hline
BUPA Liver  Disorders                       & -0.29               & \textbf{0}              & -0.88                  \\ \hline
Ozon Level  Detection                       & -0.33               & -0.06                   & \textbf{0}             \\ \hline
Parkinsons Data Set                         & -2.05               & -1.54                   & \textbf{0}             \\ \hline
Pima Indians  Diabetes                      & \textbf{0}          & -0.39                   & -0.39                  \\ \hline
Thyrod gland  data                          & \textbf{0}          & \textbf{0}              & \textbf{0}             \\ \hline
Vowel Recognition Data                      & -0.1                & -0.1                    & \textbf{0}             \\ \hline\hline
\multicolumn{1}{|c|}{\textbf{Largest Difference}}      & -4.28               & -4.28                   & \textbf{-1.38}         \\ \hline
\multicolumn{1}{|c|}{\textbf{Average Difference}}      & -0.66            & -0.46                 & \textbf{-0.27}       \\ \hline
\multicolumn{1}{|c|}{\textbf{No. of bests}} & 8                   & 8                       & \textbf{12}            \\ \hline
\end{tabular}
\end{table}

Table \ref{ClassSpeed} compares the number of linear programs solved (\textit{LPs}) and the processing time in seconds (\textit{sec}) for each algorithm, with best results in bold face. Algorithm 2E1 requires the fewest LP solutions and the least processing time in all cases. Table \ref{ClassImpSpeed} reports the percentage reduction in number of LPs and processing time for Algorithm 2E1 vs. the comparators for each dataset. For example, the processing times for Algorithms 2($\infty$) and 2E1 are 3.1 and 0.3 seconds respectively for the first dataset in Table \ref{ClassSpeed}. 
Therefore, Arithm 2E1 provides a $90.32\%$ reduction in processing time when compared to Algorithm 2($\infty$) for the first dataset in Table \ref{ClassImpSpeed}. The last two rows of Table \ref{ClassImpSpeed} show the maximum and average speed improvements obtained by Algorithm 2E1 when compared with the  other two algorithms. Algorithm 2E1 reduces the average number of LPs by about $96.56\%$ and $67.90\%$ when compared to Algorithms 2($\infty$) and 2(1), respectively, and also reduces the solution time by $94.77\%$ and $71.31\%$ when compared to  Algorithm 2($\infty$) and Algorithm 2(1) respectively.


Algorithm 2E1 thus reduces solution time by an average 70-90\% when compared to the existing algorithms, while also marginally improves the total accuracy on  average.


\begin{table}[hbt!]
\centering
\scriptsize
\caption{Algorithm speeds (number of LPs, seconds) for binary classification.}
\label{ClassSpeed}
\begin{tabular}{|l|c|c|c|c|c|c|}
\hline
\multicolumn{1}{|c|}{\multirow{2}{*}{\textbf{Dataset}}} & \multicolumn{2}{c|}{\textbf{Algorithm 2($\infty$)}} & \multicolumn{2}{c|}{\textbf{Algorithm 2(1)}} & \multicolumn{2}{c|}{\textbf{Algorithm 2E1}} \\ \cline{2-7} 
\multicolumn{1}{|c|}{}                                  & \textbf{LPs}        & \textbf{Sec}       & \textbf{LPs}          & \textbf{Sec}         & \textbf{LPs}         & \textbf{Sec}         \\ \hline
Breast Cancer Wisconsin Original                        & 100                 & 3.1                & 12                    & 0.7                  & \textbf{6}           & \textbf{0.3}         \\ \hline
Cardiotocography (Type 1 vs.others)                     & 35                  & 4.3                & 4                     & 0.6                  & \textbf{2}           & \textbf{0.3}         \\ \hline
Cardiotocography (Type 2 vs.others)                     & 67                  & 5.1                & 4                     & 0.6                  & \textbf{2}           & \textbf{0.3}         \\ \hline
Car Evaluation (Type 1 vs.others)                       & 2583                & 350.6              & 372                   & 28.1                & \textbf{8}           & \textbf{0.7}         \\ \hline
Car Evaluation (Type 2 vs.others)                       & 462                 & 28.7               & 69                    & 4.5                  & \textbf{1}           & \textbf{0.1}         \\ \hline
Car Evaluation (Type 3 vs.others)                       & 1458                & 183.2              & 203                   & 13.8                 & \textbf{10}          & \textbf{0.9}                  \\ \hline
Car Evaluation (Type 4 vs.others)                       & 205                 & 16.7               & 33                    & 2.3                  & \textbf{5}           & \textbf{0.4}         \\ \hline
Climate Model Simulation Crashes                        & 77                  & 3.2                & 7                     & 0.3                  & \textbf{5}           & \textbf{0.2}                  \\ \hline
Credit Approval                                         & 1237                & 37.1               & 87                    & 7.3                  & \textbf{1}           & \textbf{0.1}         \\ \hline
Diabetic Retinopathy  Debrecen                          & 4723                & 625.2              & 238                   & 25.6                 & \textbf{16}                   & \textbf{1.4}         \\ \hline
Glass Identification                                    & 292                 & 3.1                & 33                    & 1.5                  & \textbf{7}           & \textbf{0.3}         \\ \hline
Heart Disease Cleveland                                 & 404                 & 5.3                & 31                    & 1.0                  & \textbf{6}           & \textbf{0.2}         \\ \hline
Hepatitis Domain                                        & 16                  & 0.4                & 4                     & 0.2                  & \textbf{2}           & \textbf{0.1}         \\ \hline
Ionosphere                                              & 104                 & 4.2                & 7                     & \textbf{0.3}                  & \textbf{5}           & \textbf{0.3}         \\ \hline
BUPA Liver  Disorders                                   & 581                 & 14.8               & 83                    & 2.7                  & \textbf{10}          & \textbf{0.3}         \\ \hline
Ozon Level  Detection                                   & 542                 & 159.7              & 14                    & 4.1                  & \textbf{7}           & \textbf{2.1}         \\ \hline
Parkinsons Data Set                                     & 103                 & 2.7                & 7                     & \textbf{0.4}                 & \textbf{4}           & \textbf{0.4}         \\ \hline
Pima Indians  Diabetes                                  & 1325                & 41.7               & 153                   & 6.7                  & \textbf{12}          & \textbf{0.5}         \\ \hline
Thyrod gland  data                                      & 61                  & 1.3                & 12                    & 0.4                  & \textbf{6}           & \textbf{0.2}         \\ \hline
Vowel Recognition Data                                  & 146                 & 4.0                & 14                    & 0.5                  & \textbf{7}           & \textbf{0.3}         \\ \hline
\end{tabular}
\end{table}

\begin{table}[hbt!]
\centering
\scriptsize
\caption{Percentage reductions in LPs solved and solution times by using Algorithm 2E1 instead of original Algorithms 2($\infty$) and 2(1) for classification.}
\label{ClassImpSpeed}
\begin{tabular}{|l|c|c|c|c|}
\hline
\multicolumn{1}{|c|}{\multirow{2}{*}{\textbf{Dataset}}} & \multicolumn{2}{c|}{\textbf{\begin{tabular}[c]{@{}c@{}}Algorithm 2E1 vs. \\ Algorithm 2($\infty$)\end{tabular}}} & \multicolumn{2}{c|}{\textbf{\begin{tabular}[c]{@{}c@{}}Algorithm 2E1 vs. \\ Algorithm 2(1)\end{tabular}}} \\ \cline{2-5} 
\multicolumn{1}{|c|}{}                                  & \textbf{LPs}                                            & \textbf{Sec}                                           & \textbf{LPs}                                        & \textbf{Sec}                                        \\ \hline
Breast Cancer Wisconsin Original                        & 94                                                      & 90.32                                                  & 50                                                  & 57.14                                               \\ \hline
Cardiotocography (Type 1 vs.others)                     & 94.29                                                   & 93.02                                                  & 50                                                  & 50                                                  \\ \hline
Cardiotocography (Type 2 vs.others)                     & 97.01                                                   & 94.12                                                  & 50                                                  & 50                                                  \\ \hline
Car Evaluation (Type 1 vs.others)                       & 99.69                                                   & 99.80                                                  & 97.85                                               & 97.51                                               \\ \hline
Car Evaluation (Type 2 vs.others)                       & 99.78                                                   & 99.65                                                  & 98.55                                               & 97.78                                               \\ \hline
Car Evaluation (Type 3 vs.others)                       & 99.31                                                   & 99.51                                                  & 95.07                                               & 93.48                                               \\ \hline
Car Evaluation (Type 4 vs.others)                       & 97.56                                                   & 97.60                                                  & 84.85                                               & 82.61                                               \\ \hline
Climate Model Simulation Crashes                        & 93.51                                                   & 93.75                                                  & 28.57                                               & 33.33                                               \\ \hline
Credit Approval                                         & 99.92                                                   & 99.73                                                  & 98.85                                               & 98.63                                               \\ \hline
Diabetic Retinopathy  Debrecen                          & 99.66                                                   & 99.78                                                  & 93.28                                               & 94.53                                               \\ \hline
Glass Identification                                    & 97.60                                                   & 90.32                                                  & 78.79                                               & 80                                                  \\ \hline
Heart Disease Cleveland                                 & 98.51                                                   & 96.23                                                  & 80.64                                               & 80                                                  \\ \hline
Hepatitis Domain                                        & 87.5                                                    & 75                                                     & 50                                                  & 50                                                  \\ \hline
Ionosphere                                              & 95.19                                                   & 97.62                                                  & 28.57                                               & 66.67                                               \\ \hline
BUPA Liver  Disorders                                   & 98.28                                                   & 97.97                                                  & 87.95                                               & 88.89                                               \\ \hline
Ozon Level  Detection                                   & 98.71                                                   & 98.68                                                  & 50                                                  & 48.78                                               \\ \hline
Parkinsons Data Set                                     & 96.12                                                   & 96.30                                                  & 42.86                                               & 74.36                                               \\ \hline
Pima Indians  Diabetes                                  & 99.09                                                   & 98.80                                                  & 92.16                                               & 92.54                                               \\ \hline
Thyrod gland  data                                      & 90.16                                                   & 84.61                                                  & 50                                                  & 50                                                  \\ \hline
Vowel Recognition Data                                  & 95.20                                                   & 92.5                                                   & 50                                                  & 40                                                  \\ \hline\hline
\multicolumn{1}{|c|}{\textbf{Maximum Percent Reduction}}                  & 99.92                                                   & 99.80                                                  & 98.85                                               & 98.63                                               \\ \hline
\multicolumn{1}{|c|}{\textbf{Average Percent Reduction}}                  & 96.56                                                   & 94.77                                                  & 67.90                                               & 71.31                                               \\ \hline
\end{tabular}
\end{table}

\section{Application: Sparse Recovery in Compressive Sensing}\label{USLE}

Finding a sparse solution to an underdetermined system of linear equations in unbounded variables is an essential step in compressive sensing (CS) \cite{candes2008introduction}. A common way to compress an input signal vector $\mathbf{x}$ in CS is to multiply it by a random matrix $\mathbf{A}$ of size $m \times n$ where $m<<n$ to yield a compressed vector $\mathbf{b}$ of size $m$, which is transmitted or stored. The input vector is assumed to be sparse, so the decompression or recovery process tries to find a sparse solution given only $\mathbf{A}$ and the compressed vector $\mathbf{b}$. The CS \textit{sparse recovery problem} is to find a sparse solution $\mathbf{y}$ to the underdetermined system $\mathbf{Ay=b}$, where $\mathbf{y}$ is of size $n \times 1$ and all of its components are unbounded. 

Numerous algorithms have been developed to solve the sparse recovery problem \cite{qaisar2013compressive,meenakshi2015survey,salahdine2016survey}. Recently developed MAX FS based algorithms (Methods B, C, and M)  provide better quality solutions than the state-of-the-art algorithms in certain CS scenarios, though they can be comparatively slow \cite{firouzeh2020maximum}. 
Several choices from various families of random matrices can result in a dense $mathbf{A}$, and thus the recovery process becomes more suitable for the extended algorithms introduced here. This section examines whether the new algorithm extensions can increase the speed of the MAX FS based algorithms for the sparse recovery problem without loss of solution quality.

\subsection{MAX FS for Sparse Recovery}

Methods B, C, and M (Section \ref{findingsparsesolutions}) have been applied to CS sparse recovery and are observed to provide improvements over the state of the art in certain scenarios. Previous work \cite{firouzeh2020maximum} showed Method M to be one of the better methods for CS sparse recovery. So we combine it with our two extensions to create Method ME1E2($\ell$). The steps are summarized in Fig. \ref{ME1E2}. 

\begin{algorithm}
\captionof{figure}{Algorithm ME1E2($\ell$)}
\begin{algorithmic} [t]
\label{ME1E2}
\footnotesize
	\sffamily
\STATE $\textit{MINULR} \leftarrow \emptyset$ \\%
\text{Set up the sparse recovery LP.}\\
\textbf{do:}\\
$\hspace{1em}$\text{Solve LP.}\\
$\hspace{1em}$\text{Construct }\textit{CandidateSet}: variables $x_j$ not in \textit{MINULR} that have nonzero $\mid u_j - v_j \mid$.\\
\hspace{1em}\textbf{if} $\mid \textit{CandidateSet} \mid $ \text{$=  \phi$} \textbf{ then} \text{exit}.\\
\hspace{1em}\textbf{if } first iteration and $\mid \textit{CandidateSet} \mid $ \text{$\leq  \ell$} \textbf{ then}\\
\hspace{2em}\text{Add all candidates to} \textit{MINULR} \text{and exit}.\\
\hspace{1em}\textbf{end if}\\
\hspace{1em}Sort candidates from largest to smallest $\mid u_j - v_j \mid$.\\
\hspace{1em}Select top $p$ candidates for removal based on first abrupt change in $\mid u_j - v_j \mid$. \\%
\hspace{1em}Add the $p$ selected candidates to \textit{MINULR} and set the objective coefficients of the\\
\hspace{1em}associated $u_j,v_j$ pairs to 0.1. \\%
$\textbf{end do}$
\STATE \textbf{OUTPUT}: \textit{MINULR} is small number of nonzeros forming a support for the system of equations.  
\end{algorithmic}
\end{algorithm}

Method ME1E2($\ell$) assumes BP failure if the length of the first candidate list is greater than the defined threshold $\ell$. We set $\ell$ to $m-3$, as done for Method M in \cite{firouzeh2020maximum}.

\subsection{Comparators} 
Method ME1E2($\ell$) is compared with the original Methods B(2), C(2), and M(2). We choose list length $k=2$ to provide fast solutions with good solution quality ($k$ is typically $1-7$ \cite{chinneck2001fast}). Postprocessing is not used in any of these three methods.

\subsection{Evaluation Metric}
In sparse recovery, if the number of nonzeros in recovered signal $y$ equals the number of nonzeros in the input signal $x$, then $y$ and $x$ are usually identical. For this reason our main metric relates to the \textit{sparsity} (number of nonzeros) of $y$ vs. $x$. The \textit{critical sparsity} is the largest number of nonzeros in the input vector $x$ for which the sparsity of $y$ reliably equals the sparsity of $x$. Higher critical sparsity means that input vectors that have been more  compressed, or that they have more nonzeros and therefore can be reliably recovered. We use the critical sparsity as the evaluation metric: higher is better.

\subsection{Test Models} 
We construct examples for which the solution is known. The system is $\mathbf{A}_{128 \times 256}\mathbf{x}_{256 \times 1} = \mathbf{b}_{128 \times 1}$, providing a {compression ratio} $(CR)= (1-\frac{m}{n}) \times 100$ = 50\%. $\mathbf{A}$ is a random matrix with values ranging from -10 to 10, and hence dense. The constructed input vector $\mathbf{x}$ has various numbers of nonzeros ranging from 12 to 96, randomly distributed among its elements, and with values randomly drawn from a normal distribution with mean 0 and variance 1.  

\subsection{Hardware and Software}
See Section \ref{hwsw}.

\subsection{Results and Discussion}
Table \ref{csSPARSITY} compares Algorithm ME1E2($\ell$) with existing methods, namely, Methods B(2), C(2), and M(2). Results in the table are averages over $100$ instances at each sparsity level $S$ (actual number of nonzeros in the solution). \textit{Estimated sparsity ($T$)} is the average number of nonzeros returned by each method at each value of $S$ over $100$ instances, with the number of correct results shown in brackets. Results that are accurate in all $100$ cases are shown in boldface. The critical sparsity is the largest number shown in boldface for each method. Extended Algorithm ME1E2($\ell$) has the same critical sparsity as the other methods. 

Table \ref{csNLPs} shows the number of LPs solved and the processing time for Method ME1E2($\ell$) and the comparators. The smallest number of LPs and time in seconds are bolded. Method ME1E2($\ell$) requires the fewest  LP solutions at all values of $S$. It is significantly faster than the other algorithms, as shown in Table \ref{csIMPspeed} which reports the percentage reduction in LP solutions required and processing times for Method ME1E2($\ell$) when compared to others. The last two rows of the table show the maximum and the average percent improvement. At $S = 52$, Method ME1E2($\ell$) reduces the solution time by $96.02\%$, $92.49\%$, and $92.83\%$ compared to Methods C(2), B(2), and M(2), respectively. Method ME1E2($\ell$) does not run faster than Method M(2) when the input is very sparse (up to $S =36$), but for $S > 36$ it reduces solution time when compared to M(2) by $62.32 \%$ on average . Method ME1E2($\ell$) provides a $98.82\%$, $97.55\%$, and $66.61\%$ reduction in LPs solved when compared to Methods C(2), B(2), and M(2), respectively.It reduces the average solution time by $95.95\%$, $92.86\%$, and $62.32\%$ when compared to Methods C(2), B(2), and M(2), respectively. 

Method ME1E2($\ell$) reduces solution time significantly with no loss of solution quality (as measured by critical sparsity) for this CS sparse recovery scenario.

\begin{table}[H]
\centering
\scriptsize
\caption{Average and critical sparsities of Method ME1E2($\ell$) vs. Methods C(2), B(2) and M(2) for random matrices of size $m=128$ and $n=256$ in CS sparse recovery.}
\label{csSPARSITY}
\begin{tabular}{|c|cccc|}
\hline
\multirow{2}{*}{\textbf{S}} & \multicolumn{4}{c|}{\textbf{Estimated sparsity (T)}}                                   \\ \cline{2-5} 
                            & \textbf{Method C(2)} & \textbf{Method B(2)} & \textbf{Method M(2)} & \textbf{Method ME1E2($\ell$)} \\ \hline\hline
\textbf{12}                 & \textbf{12}$^{(100)}$        & \textbf{12}$^{(100)}$        & \textbf{12}$^{(100)}$        & \textbf{12}$^{(100)}$                 \\ \hline
\textbf{16}                 & \textbf{16}$^{(100)}$        & \textbf{16}$^{(100)}$        & \textbf{16}$^{(100)}$        & \textbf{16}$^{(100)}$                 \\ \hline
\textbf{20}                 & \textbf{20}$^{(100)}$        & \textbf{20}$^{(100)}$        & \textbf{20}$^{(100)}$        & \textbf{20}$^{(100)}$                 \\ \hline
\textbf{24}                 & \textbf{24}$^{(100)}$        & \textbf{24}$^{(100)}$        & \textbf{24}$^{(100)}$        & \textbf{24}$^{(100)}$                 \\ \hline
\textbf{28}                 & \textbf{28}$^{(100)}$        & \textbf{28}$^{(100)}$        & \textbf{28}$^{(100)}$        & \textbf{28}$^{(100)}$                 \\ \hline
\textbf{32}                 & \textbf{32}$^{(100)}$        & \textbf{32}$^{(100)}$        & \textbf{32}$^{(100)}$        & \textbf{32}$^{(100)}$                 \\ \hline
\textbf{36}                 & \textbf{36}$^{(100)}$        & \textbf{36}$^{(100)}$        & \textbf{36}$^{(100)}$        & \textbf{36}$^{(100)}$                 \\ \hline
\textbf{40}                 & \textbf{40}$^{(100)}$        & \textbf{40}$^{(100)}$        & \textbf{40}$^{(100)}$        & \textbf{40}$^{(100)}$                 \\ \hline
\textbf{44}                 & \textbf{44}$^{(100)}$        & \textbf{44}$^{(100)}$        & \textbf{44}$^{(100)}$        & \textbf{44}$^{(100)}$                 \\ \hline
\textbf{48}                 & \textbf{48}$^{(100)}$        & \textbf{48}$^{(100)}$        & \textbf{48}$^{(100)}$        & \textbf{48}$^{(100)}$                 \\ \hline
\textbf{52}                 & \textbf{52}$^{(100)}$      & \textbf{52}$^{(100)}$        & \textbf{52}$^{(100)}$        & \textbf{52}$^{(100)}$                 \\ \hline
\textbf{56}                 & 56$^{(99)}$               & 56.03$^{(98)}$               & 56.04$^{(98)}$               & 56.10$^{(98)}$                      \\ \hline
\textbf{60}                 & 64.91$^{(90)}$               & 64.24$^{(86)}$               & 65.59$^{(86)}$               & 66.43$^{(86)}$                      \\ \hline
\textbf{64}                 & 75.49$^{(77)}$                & 76.3$^{(67)}$                & 77.32$^{(67)}$               & 77.0$^{(73)}$                      \\ \hline
\textbf{68}                 & 92.21$^{(50)}$               & 94.53$^{(45)}$               & 94.48$^{(45)}$               & 92.95$^{(39)}$                      \\ \hline
\textbf{72}                 & 108.65$^{(25)}$               & 106.87$^{(25)}$              & 107.63$^{(25)}$              & 106.77$^{(19)}$                     \\ \hline 
\textbf{76}                 & 119.55$^{(9)}$               & 121.85$^{(5)}$               & 120.92$^{(5)}$               & 120.9$^{(6)}$                     \\ \hline
\textbf{80}                 & 125.41$^{(3)}$               & 126.52$^{(4)}$               & 125.56$^{(4)}$               & 125.27$^{(3)}$                      \\ \hline
\textbf{84}                 & 127.58$^{(1)}$               & 127.55$^{(0)}$               & 127.6$^{(0)}$                  & 127.14$^{(0)}$                      \\ \hline
\textbf{88}                 & 128$^{(0)}$                  & 127.98$^{(0)}$               & 128$^{(0)}$                  & 128$^{(0)}$                      \\ \hline
\textbf{92}                 & 128$^{(0)}$                  & 128$^{(0)}$                  & 128$^{(0)}$                  & 128$^{(0)}$                      \\ \hline
\textbf{96}                 & 128$^{(0)}$                  & 127.96$^{(0)}$               & 128$^{(0)}$                  & 128$^{(0)}$                      \\ \hline
\end{tabular}
\end{table}

\begin{table}[hbt!]
\centering
\scriptsize
\caption{Number of LPs and processing times for Method ME1E2($\ell$) vs. Methods C(2), B(2) and M(2) for random matrices of size $m=128$ and $n=256$ in CS sparse recovery.}
\label{csNLPs}
\begin{tabular}{|c|c|c|c|c|c|c|c|c|}
\hline
\multirow{2}{*}{\textbf{S}} & \multicolumn{2}{c|}{\textbf{Method C(2)}} & \multicolumn{2}{c|}{\textbf{Method B(2)}} & \multicolumn{2}{c|}{\textbf{Method M(2)}} & \multicolumn{2}{c|}{\textbf{Method ME1E2($\ell$)}} \\ \cline{2-9} 
                            & \textbf{LPs}       & \textbf{Sec}       & \textbf{LPs}       & \textbf{Sec}       & \textbf{LPs}       & \textbf{Sec}       & \textbf{LPs}                & \textbf{Sec}              \\ \hline
\textbf{12}                 & 48                  & 1.44                & 23                  & 0.84                & \textbf{1}          & \textbf{0.06}       & \textbf{1}                   & \textbf{0.06}              \\ \hline
\textbf{16}                 & 64                  & 2.02                & 31                  & 1.16                & \textbf{1}          & \textbf{0.05}                & \textbf{1}                   & \textbf{0.05}              \\ \hline
\textbf{20}                 & 80.06               & 3.07                & 39                  & 1.72                & \textbf{1}          & \textbf{0.06 }               & \textbf{1}                   & \textbf{0.06}              \\ \hline
\textbf{24}                 & 96.06               & 4.25                & 47                  & 2.3                 & \textbf{1}          & \textbf{0.06 }               & \textbf{1}                   & \textbf{0.06}              \\ \hline
\textbf{28}                 & 112.24              & 5.34                & 55                  & 2.98                & \textbf{1}          & \textbf{0.06}                & \textbf{1}                   & \textbf{0.06}              \\ \hline
\textbf{32}                 & 128.82              & 5.49                & 63                  & 3.02                & \textbf{1}          & \textbf{0.07 }               & \textbf{1}                   & \textbf{0.07}              \\ \hline
\textbf{36}                 & 146.14              & 6.53                & 71                  & 3.66                & \textbf{1}          & \textbf{0.07}                & \textbf{1}                   & \textbf{0.07}              \\ \hline
\textbf{40}                 & 165.84              & 7.93                & 79                  & 4.27                & 81              & 5.01               & \textbf{1.03}                & \textbf{0.08}              \\ \hline
\textbf{44}                 & 186.56              & 7.97                & 87                  & 4.27                & 99.01               & 5.23                & \textbf{1.2}                & \textbf{0.12}              \\ \hline
\textbf{48}                 & 209.46              & 8.75                & 95                  & 4.6                 & 98                  & 5.42                & \textbf{1.49}                & \textbf{0.23}              \\ \hline
\textbf{52}                 & 227.1               & 10.06               & 103                 & 5.33                & 120.22              & 5.58                & \textbf{2.35}                & \textbf{0.40}              \\ \hline
\textbf{56}                 & 246.87              & 11.33               & 111.12              & 5.99                & 132.01              & 6.02                & \textbf{3.01}               & \textbf{0.53}              \\ \hline
\textbf{60}                 & 286.09              & 13.41               & 127.68              & 7.02                & 133.45              & 7.26                & \textbf{3.82}               & \textbf{0.73}              \\ \hline
\textbf{64}                 & 328.65              & 15.72               & 152.26              & 8.56                & 165.78              & 9.01                & \textbf{4.63}               & \textbf{0.90}              \\ \hline
\textbf{68}                 & 395.37              & 18.98               & 189.39              & 10.72               & 192.08              & 11.54               & \textbf{5.46}               & \textbf{1.12}              \\ \hline
\textbf{72}                 & 456.75              & 21.87               & 214.48              & 12.15               & 253.66              & 12.01               & \textbf{5.95}               & \textbf{1.31}              \\ \hline
\textbf{76}                 & 493                 & 23.7                & 245.06              & 13.91               & 256.08              & 14.05               & \textbf{6.36}                & \textbf{1.42}              \\ \hline
\textbf{80}                 & 505.55              & 24.02               & 254.84              & 14.37               & 259.06              & 14.53               & \textbf{6.18}               & \textbf{1.41}              \\ \hline
\textbf{84}                 & 509.31              & 24.25               & 257.09              & 14.5                & 259.24              & 14.54               & \textbf{6.85}               & \textbf{1.58}               \\ \hline
\textbf{88}                 & 511.12              & 24.28               & 258.2               & 14.55               & 260.02              & 14.68               & \textbf{6.75}               & \textbf{1.57}              \\ \hline
\textbf{92}                 & 511.15              & 24.9                & 258.18              & 14.83               & 261.01              & 15.02               & \textbf{6.87}               & \textbf{1.67}              \\ \hline
\textbf{96}                 & 511.1616            & 24.52               & 257.9394            & 14.64               & 261                 & 15.11               & \textbf{7.3}               & \textbf{1.79}              \\ \hline
\end{tabular}
\end{table}
\begin{table}[hbt!]
\centering
\scriptsize
\caption{Percent reduction in number of LPs and processing times for Method ME1E2($\ell$) vs. Methods C(2), B(2) and M(2) for random matrices of size $m=128$ and $n=256$ in CS sparse recovery.}
\label{csIMPspeed}
\begin{tabular}{|c|c|c|c|c|c|c|}
\hline
\multirow{2}{*}{\textbf{S}} & \multicolumn{2}{c|}{\textbf{\begin{tabular}[c]{@{}c@{}}Method ME1E2($\ell$) vs.\\ Method C(2)\end{tabular}}} & \multicolumn{2}{c|}{\textbf{\begin{tabular}[c]{@{}c@{}}Method ME1E2($\ell$) vs.\\ Method B(2)\end{tabular}}} & \multicolumn{2}{c|}{\textbf{\begin{tabular}[c]{@{}c@{}}Method ME1E2($\ell$) vs.\\ Method M(2)\end{tabular}}} \\ \cline{2-7} 
                            & \textbf{LPs}                                          & \textbf{Sec}                                         & \textbf{LPs}                                          & \textbf{Sec}                                         & \textbf{LPs}                                          & \textbf{Sec}                                         \\ \hline
\textbf{12}                 & 97.92                                                 & 95.83333                                             & 95.65                                                 & 92.86                                                & 0                                                     & 0                                                    \\ \hline
\textbf{16}                 & 98.44                                                 & 97.53                                                & 96.77                                                 & 95.69                                                & 0                                                     & 0                                                    \\ \hline
\textbf{20}                 & 98.75                                                 & 98.05                                                & 97.44                                                 & 96.51                                                & 0                                                     & 0                                                    \\ \hline
\textbf{24}                 & 98.96                                                 & 98.59                                                & 97.87                                                 & 97.39                                                & 0                                                     & 0                                                    \\ \hline
\textbf{28}                 & 99.11                                                 & 98.88                                                & 98.18                                                 & 97.99                                                & 0                                                     & 0                                                    \\ \hline
\textbf{32}                 & 99.22                                                 & 98.73                                                & 98.41                                                 & 97.68                                                & 0                                                     & 0                                                    \\ \hline
\textbf{36}                 & 99.32                                                 & 98.93                                                & 98.59                                                 & 98.09                                                & 0                                                     & 0                                                    \\ \hline
\textbf{40}                 & 99.38                                                 & 98.99                                                & 98.70                                                 & 98.13                                                & 98.73                                                 & 98.403                                               \\ \hline
\textbf{44}                 & 99.36                                                 & 98.49                                                & 98.62                                                 & 97.19                                                & 98.79                                                 & 97.70                                                \\ \hline
\textbf{48}                 & 99.29                                                 & 97.37                                                & 98.43                                                 & 95                                                   & 98.48                                                 & 95.76                                                \\ \hline
\textbf{52}                 & 98.96                                                 & 96.02                                                & 97.72                                                 & 92.49                                                & 98.04                                                 & 92.83                                                \\ \hline
\textbf{56}                 & 98.78                                                 & 95.32                                                & 97.29                                                 & 91.15                                                & 97.72                                                 & 91.20                                                \\ \hline
\textbf{60}                 & 98.66                                                 & 94.56                                                & 97.01                                                 & 89.60                                                & 97.14                                                 & 89.94                                                \\ \hline
\textbf{64}                 & 98.59                                                 & 94.28                                                & 96.96                                                 & 89.49                                                & 97.21                                                 & 90.01                                                \\ \hline
\textbf{68}                 & 98.62                                                 & 94.10                                                & 97.12                                                 & 89.55                                                & 97.16                                                 & 90.29                                                \\ \hline
\textbf{72}                 & 98.69                                                 & 94.01                                                & 97.23                                                 & 89.22                                                & 97.65                                                 & 89.09                                                \\ \hline
\textbf{76}                 & 98.71                                                 & 94.01                                                & 97.40                                                 & 89.79                                                & 97.52                                                 & 89.89                                                \\ \hline
\textbf{80}                 & 98.78                                                 & 94.13                                                & 97.57                                                 & 90.19                                                & 97.61                                                 & 90.29                                                \\ \hline
\textbf{84}                 & 98.65                                                 & 93.48                                                & 97.34                                                 & 89.10                                                & 97.36                                                 & 89.13                                                \\ \hline
\textbf{88}                 & 98.68                                                 & 93.53                                                & 97.39                                                 & 89.21                                                & 97.40                                                 & 89.30                                                \\ \hline
\textbf{92}                 & 98.66                                                 & 93.41                                                & 97.34                                                 & 88.94                                                & 97.37                                                 & 89.08                                                \\ \hline
\textbf{96}                 & 98.577                                                & 92.670                                               & 97.17                                                 & 87.77                                                & 97.20                                                 & 88.15                                                \\ \hline \hline
\textbf{Maximum}            & 99.38                                                 & 98.99                                                & 98.70                                                 & 98.13                                                & 98.79                                                 & 98.40                                                \\ \hline
\textbf{Average}            & 98.82                                                 & 95.95                                                & 97.55                                                 & 92.86                                                & 66.61                                                 & 62.32                                                \\ \hline
\end{tabular}
\end{table}

\newpage
\section{Conclusions} \label{Conc}

We propose two new extensions to Chinneck's original MAX FS solution algorithms \cite{chinneck1996effective} \cite{chinneck2001fast} for use with dense constraint matrices. These extensions increase solution speed greatly with no loss of solution quality when tested in two applications: classification, and sparse recovery in compressive sensing:
\begin{itemize}
\item \textit{Classification}: The new Algorithm 2E1 reduces the mean solution time by about $94.77\%$ and $71.31\%$ when compared to Algorithms 2($\infty$) and 2(1), respectively, while slightly improving accuracy on average. The better classification accuracies reflect larger feasible subsets found by the algorithm over its predecessors.
\item \textit{Sparse recovery in compressive sensing}: The new Method ME1E2($\ell$) reduces processing time by $95.95\%$, $92.86\%$ and $62.32\%$ when compared with Methods C, B and M, on average, with no reduction in the critical sparsity.
\end{itemize}

We expect that these results will generalize to other applications that have dense constraint matrices. We plan to investigate the use of the extended algorithms in the recovery of compressively sensed biomedical signals and dictionary learning. We also intend to extend the basic ideas to less dense constraint matrices, including those found in general LPs.





\bibliographystyle{elsarticle-num-names}
\bibliography{Manuscript.bib}

\begin{thebibliography}{34}
\expandafter\ifx\csname natexlab\endcsname\relax\def\natexlab#1{#1}\fi
\providecommand{\url}[1]{\texttt{#1}}
\providecommand{\href}[2]{#2}
\providecommand{\path}[1]{#1}
\providecommand{\DOIprefix}{doi:}
\providecommand{\ArXivprefix}{arXiv:}
\providecommand{\URLprefix}{URL: }
\providecommand{\Pubmedprefix}{pmid:}
\providecommand{\doi}[1]{\href{http://dx.doi.org/#1}{\path{#1}}}
\providecommand{\Pubmed}[1]{\href{pmid:#1}{\path{#1}}}
\providecommand{\bibinfo}[2]{#2}
\ifx\xfnm\relax \def\xfnm[#1]{\unskip,\space#1}\fi
\bibitem[{Amaldi et~al.(1999)Amaldi, Pfetsch, and Trotter}]{One}
\bibinfo{author}{E.~Amaldi}, \bibinfo{author}{M.~E. Pfetsch},
  \bibinfo{author}{L.~E. Trotter},
\newblock \bibinfo{title}{Some structural and algorithmic properties of the
  maximum feasible subsystem problem},
\newblock \bibinfo{journal}{In International Conference on Integer Programming
  and Combinatorial Optimization}  (\bibinfo{year}{1999})
  \bibinfo{pages}{45–59}.
\bibitem[{Amaldi(1994)}]{amaldi1994finding}
\bibinfo{author}{E.~Amaldi}, \bibinfo{title}{From finding maximum feasible
  subsystems of linear systems to feedforward neural network design},
  \bibinfo{type}{{PhD} dissertation}, Swiss Federal Institute of Technology at
  Lausanne (EPFL), \bibinfo{year}{1994}.
\bibitem[{Chinneck(1996)}]{chinneck1996effective}
\bibinfo{author}{J.~W. Chinneck},
\newblock \bibinfo{title}{An effective polynomial-time heuristic for the
  minimum-cardinality iis set-covering problem},
\newblock \bibinfo{journal}{Annals of Mathematics and Artificial Intelligence}
  \bibinfo{volume}{17} (\bibinfo{year}{1996}) \bibinfo{pages}{127--144}.
\bibitem[{Chakravarti(1994)}]{Seven}
\bibinfo{author}{N.~Chakravarti},
\newblock \bibinfo{title}{Some results concerning post-infeasibility analysis},
\newblock \bibinfo{journal}{Oper. Res.} \bibinfo{volume}{73}
  (\bibinfo{year}{1994}) \bibinfo{pages}{139–143}.
\bibitem[{Amaldi and Kann(1994)}]{Eight}
\bibinfo{author}{E.~Amaldi}, \bibinfo{author}{V.~Kann},
\newblock \bibinfo{title}{The complexity and approximability of finding maximum
  feasible subsystems of linear relations},
\newblock \bibinfo{journal}{Oper. Res.} \bibinfo{volume}{147}
  (\bibinfo{year}{1994}) \bibinfo{pages}{181–210}.
\bibitem[{Sankaran(1993)}]{Nine}
\bibinfo{author}{J.~K. Sankaran},
\newblock \bibinfo{title}{A note on resolving infeasibility in linear programs
  by constraint relaxation},
\newblock \bibinfo{journal}{Oper. Res.} \bibinfo{volume}{13}
  (\bibinfo{year}{1993}) \bibinfo{pages}{19--20}.
\bibitem[{Bennett and Bredensteiner(1997)}]{Three}
\bibinfo{author}{K.~P. Bennett}, \bibinfo{author}{E.~J. Bredensteiner},
\newblock \bibinfo{title}{A parametric optimization method for machine
  learning},
\newblock \bibinfo{journal}{INFORMS Journal on Computing} \bibinfo{volume}{9}
  (\bibinfo{year}{1997}) \bibinfo{pages}{311--318}.
\bibitem[{Mangasarian(1994)}]{Four}
\bibinfo{author}{O.~L. Mangasarian},
\newblock \bibinfo{title}{Misclassification minimization},
\newblock \bibinfo{journal}{Journal of Global Optimization} \bibinfo{volume}{5}
  (\bibinfo{year}{1994}) \bibinfo{pages}{309--323}.
\bibitem[{Rossi et~al.(2001)Rossi, Smriglio, and Sassano}]{Five}
\bibinfo{author}{F.~Rossi}, \bibinfo{author}{S.~Smriglio},
  \bibinfo{author}{A.~Sassano},
\newblock \bibinfo{title}{Models and algorithms for terrestrial digital},
\newblock \bibinfo{journal}{Annals of Operations Research}
  \bibinfo{volume}{107} (\bibinfo{year}{2001}) \bibinfo{pages}{267–283}.
\bibitem[{Wagner and Elber(2004)}]{Six}
\bibinfo{author}{M.~Wagner}, \bibinfo{author}{R.~Elber},
\newblock \bibinfo{title}{Large-scale linear programming techniques for the
  design},
\newblock \bibinfo{journal}{Annals of Operations Research}
  \bibinfo{volume}{318} (\bibinfo{year}{2004}) \bibinfo{pages}{301–318}.
\bibitem[{Firouzeh et~al.(2020)Firouzeh, Chinneck, and
  Rajan}]{firouzeh2020maximum}
\bibinfo{author}{F.~F. Firouzeh}, \bibinfo{author}{J.~W. Chinneck},
  \bibinfo{author}{S.~Rajan},
\newblock \bibinfo{title}{Maximum feasible subsystem algorithms for recovery of
  compressively sensed speech},
\newblock \bibinfo{journal}{IEEE Access} \bibinfo{volume}{8}
  (\bibinfo{year}{2020}) \bibinfo{pages}{82539--82550}.
\bibitem[{Chinneck(2007)}]{chinneck2008feasibility}
\bibinfo{author}{J.~W. Chinneck}, \bibinfo{title}{Feasibility and Infeasibility
  in Optimization: Algorithms and Computational Methods}, volume
  \bibinfo{volume}{118}, \bibinfo{publisher}{Springer Science \& Business
  Media}, \bibinfo{year}{2007}.
\bibitem[{Amaldi(2003)}]{amaldi2003maximum}
\bibinfo{author}{E.~Amaldi},
\newblock \bibinfo{title}{The maximum feasible subsystem problem and some
  applications},
\newblock \bibinfo{journal}{Modelli e Algoritmi per l’Ottimizzazione di
  Sistemi Complessi, A. Agnetis and G. D. Pillo, eds., Pitagora Editrice,
  Bologna, Italy}  (\bibinfo{year}{2003}) \bibinfo{pages}{31–69}.
\bibitem[{Mangasarian(1994)}]{mangasarian1994misclassification}
\bibinfo{author}{O.~L. Mangasarian},
\newblock \bibinfo{title}{Misclassification minimization},
\newblock \bibinfo{journal}{Journal of Global Optimization} \bibinfo{volume}{5}
  (\bibinfo{year}{1994}) \bibinfo{pages}{309--323}.
\bibitem[{Mangasarian(1996)}]{mangasarian1996machine}
\bibinfo{author}{O.~Mangasarian},
\newblock \bibinfo{title}{Machine learning via polyhedral concave
  minimization},
\newblock in: \bibinfo{booktitle}{Applied Mathematics and Parallel Computing},
  \bibinfo{publisher}{Springer}, \bibinfo{year}{1996}, pp.
  \bibinfo{pages}{175--188}.
\bibitem[{Bennett and Bredensteiner(1997)}]{bennett1997parametric}
\bibinfo{author}{K.~P. Bennett}, \bibinfo{author}{E.~J. Bredensteiner},
\newblock \bibinfo{title}{A parametric optimization method for machine
  learning},
\newblock \bibinfo{journal}{INFORMS Journal on Computing} \bibinfo{volume}{9}
  (\bibinfo{year}{1997}) \bibinfo{pages}{311--318}.
\bibitem[{Parker(1995)}]{parker1995set}
\bibinfo{author}{M.~R. Parker}, \bibinfo{title}{A set covering approach to
  infeasibility analysis of linear programming problems and related issues},
  Ph.D. thesis, University of Colorado at Denver Denver, Colorado,
  \bibinfo{year}{1995}.
\bibitem[{Parker and Ryan(1996)}]{parker1996finding}
\bibinfo{author}{M.~Parker}, \bibinfo{author}{J.~Ryan},
\newblock \bibinfo{title}{Finding the minimum weight iis cover of an infeasible
  system of linear inequalities},
\newblock \bibinfo{journal}{Annals of Mathematics and Artificial Intelligence}
  \bibinfo{volume}{17} (\bibinfo{year}{1996}) \bibinfo{pages}{107--126}.
\bibitem[{Pfetsch(2003)}]{pfetsch2003maximum}
\bibinfo{author}{M.~E. Pfetsch}, \bibinfo{title}{The maximum feasible subsystem
  problem and vertex-facet incidences of polyhedra}, Ph.D. thesis, TU Berlin,
  Berlin, \bibinfo{year}{2003}.
\bibitem[{Pfetsch(2008)}]{pfetsch2008branch}
\bibinfo{author}{M.~E. Pfetsch},
\newblock \bibinfo{title}{Branch-and-cut for the maximum feasible subsystem
  problem},
\newblock \bibinfo{journal}{SIAM Journal on Optimization} \bibinfo{volume}{19}
  (\bibinfo{year}{2008}) \bibinfo{pages}{21--38}.
\bibitem[{Chinneck(2001)}]{chinneck2001fast}
\bibinfo{author}{J.~W. Chinneck},
\newblock \bibinfo{title}{Fast heuristics for the maximum feasible subsystem
  problem},
\newblock \bibinfo{journal}{INFORMS Journal on Computing} \bibinfo{volume}{13}
  (\bibinfo{year}{2001}) \bibinfo{pages}{210--223}.
\bibitem[{Jokar and Pfetsch(2008)}]{jokar2008exact}
\bibinfo{author}{S.~Jokar}, \bibinfo{author}{M.~E. Pfetsch},
\newblock \bibinfo{title}{Exact and approximate sparse solutions of
  underdetermined linear equations},
\newblock \bibinfo{journal}{SIAM Journal on Scientific Computing}
  \bibinfo{volume}{31} (\bibinfo{year}{2008}) \bibinfo{pages}{23--44}.
\bibitem[{Chen et~al.(2001)Chen, Donoho, and Saunders}]{chen2001atomic}
\bibinfo{author}{S.~S. Chen}, \bibinfo{author}{D.~L. Donoho},
  \bibinfo{author}{M.~A. Saunders},
\newblock \bibinfo{title}{Atomic decomposition by basis pursuit},
\newblock \bibinfo{journal}{SIAM review} \bibinfo{volume}{43}
  (\bibinfo{year}{2001}) \bibinfo{pages}{129--159}.
\bibitem[{Chinneck(2018)}]{ChinneckMAXFS}
\bibinfo{author}{J.~W. Chinneck},
\newblock \bibinfo{title}{Sparse solutions of linear systems via maximum
  feasible subsets},
\newblock \bibinfo{journal}{Les Cahiers du GERAD G-2018-104}
  (\bibinfo{year}{2018}).
\bibitem[{Killick et~al.(2012)Killick, Fearnhead, and
  Eckley}]{killick2012optimal}
\bibinfo{author}{R.~Killick}, \bibinfo{author}{P.~Fearnhead},
  \bibinfo{author}{I.~A. Eckley},
\newblock \bibinfo{title}{Optimal detection of changepoints with a linear
  computational cost},
\newblock \bibinfo{journal}{Journal of the American Statistical Association}
  \bibinfo{volume}{107} (\bibinfo{year}{2012}) \bibinfo{pages}{1590--1598}.
\bibitem[{Chinneck(2009)}]{chinneck2009tailoring}
\bibinfo{author}{J.~W. Chinneck},
\newblock \bibinfo{title}{Tailoring classifier hyperplanes to general metrics},
\newblock in: \bibinfo{booktitle}{Operations research and
  cyber-infrastructure}, \bibinfo{publisher}{Springer}, \bibinfo{year}{2009},
  pp. \bibinfo{pages}{365--387}.
\bibitem[{Chinneck(2012)}]{chinneck2012integrated}
\bibinfo{author}{J.~W. Chinneck},
\newblock \bibinfo{title}{Integrated classifier hyperplane placement and
  feature selection},
\newblock \bibinfo{journal}{Expert Systems with Applications}
  \bibinfo{volume}{39} (\bibinfo{year}{2012}) \bibinfo{pages}{8193--8203}.
\bibitem[{Silva(2017)}]{silva2017optimization}
\bibinfo{author}{A.~P.~D. Silva},
\newblock \bibinfo{title}{Optimization approaches to supervised
  classification},
\newblock \bibinfo{journal}{European Journal of Operational Research}
  \bibinfo{volume}{261} (\bibinfo{year}{2017}) \bibinfo{pages}{772--788}.
\bibitem[{Dua and Graff(2017)}]{Dua:2019}
\bibinfo{author}{D.~Dua}, \bibinfo{author}{C.~Graff}, \bibinfo{title}{{UCI}
  machine learning repository}, \bibinfo{year}{2017}. \URLprefix
  \url{http://archive.ics.uci.edu/ml}.
\bibitem[{{Mosek ApS }(2018)}]{MOSEK}
\bibinfo{author}{{Mosek ApS }}, \bibinfo{title}{The mosek optimization
  software}, \bibinfo{year}{2018}. \URLprefix \url{https://www.mosek.com},
  \bibinfo{note}{last accessed 2021}.
\bibitem[{Cand{\`e}s and Wakin(2008)}]{candes2008introduction}
\bibinfo{author}{E.~J. Cand{\`e}s}, \bibinfo{author}{M.~B. Wakin},
\newblock \bibinfo{title}{An introduction to compressive sampling},
\newblock \bibinfo{journal}{IEEE signal processing magazine}
  \bibinfo{volume}{25} (\bibinfo{year}{2008}) \bibinfo{pages}{21--30}.
\bibitem[{Qaisar et~al.(2013)Qaisar, Bilal, Iqbal, Naureen, and
  Lee}]{qaisar2013compressive}
\bibinfo{author}{S.~Qaisar}, \bibinfo{author}{R.~M. Bilal},
  \bibinfo{author}{W.~Iqbal}, \bibinfo{author}{M.~Naureen},
  \bibinfo{author}{S.~Lee},
\newblock \bibinfo{title}{Compressive sensing: From theory to applications, a
  survey},
\newblock \bibinfo{journal}{Journal of Communications and networks}
  \bibinfo{volume}{15} (\bibinfo{year}{2013}) \bibinfo{pages}{443--456}.
\bibitem[{Meenakshi(2015)}]{meenakshi2015survey}
\bibinfo{author}{S.~B. Meenakshi},
\newblock \bibinfo{title}{A survey of compressive sensing based greedy pursuit
  reconstruction algorithms},
\newblock \bibinfo{journal}{International Journal of Image, Graphics and Signal
  Processing} \bibinfo{volume}{7} (\bibinfo{year}{2015})
  \bibinfo{pages}{1--10}.
\bibitem[{Salahdine et~al.(2016)Salahdine, Kaabouch, and
  El~Ghazi}]{salahdine2016survey}
\bibinfo{author}{F.~Salahdine}, \bibinfo{author}{N.~Kaabouch},
  \bibinfo{author}{H.~El~Ghazi},
\newblock \bibinfo{title}{A survey on compressive sensing techniques for
  cognitive radio networks},
\newblock \bibinfo{journal}{Physical Communication} \bibinfo{volume}{20}
  (\bibinfo{year}{2016}) \bibinfo{pages}{61--73}.

\end{thebibliography}








\end{document}